\documentclass[twoside,notitlepage,12pt]{article}

\usepackage{amssymb}
\usepackage[leqno]{amsmath}
\usepackage{amsfonts}
\usepackage{amsopn}
\usepackage{amstext}
\usepackage{amsthm}

\usepackage[all]{xy}
\newdir{ >}{{}*!/-9pt/@{>}}

\usepackage[colorlinks]{hyperref}

\providecommand{\cal}{\mathcal}
\renewcommand{\Bbb}{\mathbb}

\newenvironment{pf}{\begin{proof}}{\end{proof}}

\newcommand{\Aaa}{{\cal{A}}}
\newcommand{\Bee}{{\cal{B}}}

\newcommand{\Dee}{{\cal{D}}}
\newcommand{\Ef}{{\cal{F}}}
\newcommand{\Gee}{{\cal{G}}}

\newcommand{\Kay}{{\cal{K}}}

\newcommand{\Nat}{{\Bbb{N}}}
\newcommand{\Qyu}{{\Bbb{Q}}}
\newcommand{\Err}{{\Bbb{R}}}

\newcommand{\lam}{{\lambda}}
\newcommand{\al}{\alpha}

\newcommand{\eps}{\varepsilon}
\renewcommand{\phi}{\varphi}
\renewcommand{\rho}{\varrho}

\newcommand{\G}{{\mathbb G}}
\newcommand{\V}{{\mathbb V}}
\newcommand{\W}{{\mathbb W}}

\newcommand{\rest}{\restriction}

\newcommand{\ntr}{{n\in\omega}}
\newcommand{\Ntr}{n\in{\Bbb{N}}}
\newcommand{\loe}{\leqslant}
\newcommand{\goe}{\geqslant}

\newcommand{\subs}{\subseteq}
\newcommand{\sups}{\supseteq}

\newcommand{\cl}{\operatorname{cl}}

\newcommand{\dens}[1]{\operatorname{dens}\left(#1\right)}

\newcommand{\id}[1]{{\operatorname{i\!d}_{#1}}} 
\newcommand{\cf}{\operatorname{cf}}

\newcommand{\poset}{{\Bbb{P}}}

\newcommand{\by}{/}

\newtheorem{tw}{Theorem}[section]
\newtheorem{wn}[tw]{Corollary}
\newtheorem{lm}[tw]{Lemma}
\newtheorem{prop}[tw]{Proposition}
\newtheorem{claim}[tw]{Claim}
\theoremstyle{definition}
\newtheorem{df}[tw]{Definition}

\newtheorem{question}[tw]{Question}
\newtheorem{problem}[tw]{Problem}
\theoremstyle{remark}

\newcommand{\setof}[2]{\{#1\colon #2\}}
\newcommand{\bigsetof}[2]{\Bigl\{#1\colon #2\Bigr\}}

\newcommand{\seq}[1]{\langle #1 \rangle}

\newcommand{\sett}[2]{\{#1\}_{#2}}
\newcommand{\sn}[1]{\{#1\}} 
\newcommand{\pair}[2]{\langle #1, #2 \rangle} 
\newcommand{\map}[3]{#1\colon #2 \to #3} 
\newcommand{\img}[2]{#1[#2]} 

\newcommand{\im}{\operatorname{im}}

\newcommand{\fra}{Fra\"iss\'e}
\newcommand{\jon}{J\'onsson}
\newcommand{\frajon}{\fra-\jon}

\providecommand{\nat}{\omega}

\newcommand{\ciag}[1]{{\sett{{#1}_n}{\ntr}}}

\newcommand{\anorm}{\|\cdot\|}
\newcommand{\norm}[1]{\|#1\|}
\newcommand{\bnorm}[1]{\Bigl\|#1\Bigr\|}

\newcommand{\abs}[1]{|#1|}

\newcommand{\usphere}[1]{\operatorname{S}_{#1}}

\newcommand{\fK}{{\mathfrak{K}}}

\newcommand{\fD}{{\mathfrak{D}}}

\newcommand{\cmp}{\circ} 

\newcommand{\cont}{\ensuremath{\mathfrak c}}

\newcommand{\til}[1]{\tilde{#1}}

\newcommand{\separator}{\begin{center}***\end{center}}

\newcommand{\ela}{{\ell_\infty}}

\newcommand{\ud}{\operatorname{UD}}
\newcommand{\udsep}{\operatorname{UD}(\mbox{sep})}
\newcommand{\udfd}{\operatorname{UD}(\mbox{fin-dim})}

\newcommand{\Gurarii}{Gurari\u\i}

\newcommand{\pogfd}{\mathfrak P\mathfrak O_{fd}}

\title{Remarks on \Gurarii\ spaces}
\author{
{Joanna Garbuli\'nska}\footnote{
Research of the first author is supported by the ESF Human Capital Operational Program grant 6/1/8.2.1./POKL/2009.
}\\
{\small Institute of Mathematics,}
{\small Jan Kochanowski University in Kielce, Poland}\\
\and
{Wies{\l}aw Kubi\'s}\footnote{
Research of the second author is supported in part
by the Grant IAA 100 190 901 and by the Institutional Research
Plan of the Academy of Sciences of Czech Republic No. AVOZ 101 905 03.
}\\
{\small Institute of Mathematics, Academy of Sciences of the Czech Republic}\\
{\small\texttt{kubis@math.cas.cz}}
}

\begin{document}

\maketitle

\begin{abstract}
We present selected known results and some of their improvements, involving \Gurarii\ spaces.
A Banach space is \emph{\Gurarii} if it has certain natural extension property for almost isometric embeddings of finite-dimensional spaces.
Deleting the word ``almost", we get the notion of a \emph{strong \Gurarii} space.
There exists a unique (up to isometry) separable \Gurarii\ space, however strong \Gurarii\ spaces cannot be separable.
The structure of the class of non-separable \Gurarii\ spaces seems to be not very well understood. We discuss some of their properties and state some open questions.
In particular, we characterize non-separable \Gurarii\ spaces in terms of skeletons of separable subspaces, we construct a non-separable \Gurarii\ space with a projectional resolution of the identity and we show that no strong \Gurarii\ space can be weakly Lindel\"of determined.
\\


\end{abstract}

\tableofcontents

\section*{Introduction}

The \emph{\Gurarii\ space}, constructed by \Gurarii~\cite{Gurarii} in 1965, is the unique separable Banach space $\G$ satisfying the following condition: Given finite-dimensional Banach spaces $X\subs Y$, given $\eps > 0$, given an isometric linear embedding $\map fX\G$ there exists an injective linear operator $\map gY\G$ extending $f$ and satisfying $\norm g \cdot \norm {g^{-1}} < 1 + \eps$.
Almost straight from this definition, it is not hard to prove that such a space is unique up to isomorphism of norm arbitrarily close to one.
Surprisingly, it has been unknown for some time whether the \Gurarii\ space is unique up to isometry; it was answered affirmatively by Lusky~\cite{Lusky} in 1976.
His proof used the method of \emph{representing matrices}, explored earlier by Lazar and Lindenstrauss~\cite{LazLin}. Very recently, Solecki and the second author~\cite{KubSol} have found a simple and elementary proof of the uniqueness of the \Gurarii\ space. We sketch the arguments in Section~\ref{Segurademolinaaa} below.

The defining condition of a \Gurarii\ space can clearly be applied to non-separable spaces, obtaining the notion of a \emph{\Gurarii\ space}. Removing $\eps$ from the definition, one gets the notion of a \emph{strong \Gurarii\ space}.
Besides their existence, not much is known about the structure of strong \Gurarii\ spaces.
Few years ago, the second author found, assuming the continuum hypothesis, a unique Banach space $V$ of density continuum and satisfying the following stronger property: every isometric embedding $\map fSV$ from a subspace of an arbitrary fixed separable space $T$ can be extended to an isometric  embedding $\map gTV$.
In fact, this is a special case of a general theory of \emph{\frajon\ limits}.
Recently, the authors of \cite{ACCGMud} developed the idea of ``generating" Banach spaces by using pushouts, finding strong \Gurarii\ spaces of arbitrarily large density above the continuum.

In this note we survey the basic properties of the separable \Gurarii\ space, we explain the pushout constructions, and we characterize \Gurarii\ spaces in terms of skeletons of separable spaces.
We also show that Banach spaces constructed by pushout iterations from finite-dimensional spaces are not universal for spaces of density $\aleph_1$. More specifically, we show that every copy of $c_0$ is complemented in such spaces.
Finally, we state some questions regarding the structure of \Gurarii\ spaces.

The paper is organized as follows.
Section~\ref{SecPrelimsss} contains the basic definitions and an overview of the Pushout Lemma, crucial for the existence of \Gurarii\ spaces.
Section~\ref{Segurademolinaaa} has a survey character. We introduce \Gurarii\ spaces, describe two natural constructions, and sketch the proof of their isometric uniqueness.
We also provide a proof of the result of Wojtaszczyk~\cite{Wojt} on 1-complemented subspaces of the \Gurarii\ space.
Section~\ref{Secnonsepeirge} studies non-separable \Gurarii\ spaces. We characterize them in terms of skeletons of separable subspaces.
As an application, we observe that no \Gurarii\ space is complemented in a $C(K)$ space and we prove that every Banach space embeds isometrically into a \Gurarii\ space of the same density.
We also show that there exists a \Gurarii\ space of density $\aleph_1$ and with a projectional resolution of the identity.
Section~\ref{SecUDudud} deals with a natural generalization of the notion of a strong \Gurarii\ space, when the class of finite-dimensional spaces is replaced by a larger class $\Kay$. The property is then called ``universal disposition for $\Kay$".
We review the ``pushout construction" which is the main tool in~\cite{ACCGMud} for constructing spaces of universal disposition for various classes.
Section~\ref{Sergowrgwrg} addresses the structure of strong \Gurarii\ spaces.
Using the fact that the \Gurarii\ space is not 1-injective for finite-dimen\-sio\-nal spaces, we observe that strong \Gurarii\ spaces cannot contain skeletons of 1-complemented separable subspaces; in particular no weakly compactly generated space can be a strong \Gurarii\ space.
We finally show, using some arguments from \cite{AviBre} that strong \Gurarii\ spaces constructed by pushout iterations in \cite{ACCGMud} have the property that every copy of $c_0$ is complemented.
Section~\ref{SecTheLastOne} contains some concluding remarks and some open questions.

\section{Preliminaries}\label{SecPrelimsss}

We shall use standard notions concerning Banach spaces and linear operators (all linear operator are, by default, bounded).
We shall consider real Banach spaces, although the result are valid for the complex case, without any significant changes.

The following well-known notion will be used throughout this paper.
Let $X$, $Y$ be Banach spaces, $\eps>0$. A linear operator $\map fXY$ is an \emph{$\eps$-isometric embedding} if
$$ (1+\eps)^{-1} \cdot \norm {x} \loe \norm {f(x)} \loe (1+\eps) \cdot \norm {x}.$$
holds for every $x\in X\setminus \sn0$.
If the above condition holds with strict inequalities, we shall say that $f$ is a \emph{strict $\eps$-isometric embedding}.
An operator $f$ is an \emph{isometric embedding} iff it is an $\eps$-isometric embedding with $\eps = 0$.
A bijective ($\eps$-)isometric embedding is called an \emph{($\eps$-)isometry}.
(The word ``isometry" always means ``linear isometry".)
Two Banach spaces are \emph{linearly isometric} if there exists a linear isometry between.
Two Banach spaces are \emph{almost linearly isometric} if for every $\eps > 0$ there exists a linear $\eps$-isometry between them.
Two norms on the same Banach space are \emph{$\eps$-equivalent} if the identity is an $\eps$-isometry.

We shall need the following simple and standard fact on extending equivalent norms.

\begin{lm}\label{Lertnort}
Let $E\subs F$ be Banach spaces, $\eps > 0$ and let $\abs{\cdot}_E$ be a norm on $E$ that is $\eps$-equivalent to the original norm of $E$ (inherited from $F$).
Then there exists a norm $\abs{\cdot}_F$ that extends $\abs{\cdot}_E$ and is $\eps$-equivalent to the original norm of $F$.
\end{lm}

\begin{pf}
Let $\anorm$ be the original norm of $F$ and let $S = \setof{\phi\in E^*}{\abs{\phi}_E=1}$ be the dual sphere in $E^*$ with respect to $\abs{\cdot}_E$.
Then $\norm\phi \loe 1+\eps$ for every $\phi\in S$.
Given $y\in F$, define
$$\abs{y}_F = \sup\setof{\psi(y)}{\psi\rest E \in S \text{ and } \norm\psi \loe 1+\eps}.$$
It is clear that $\abs{\cdot}_F$ extends $\abs{\cdot}_E$ and is $\eps$-equivalent $\anorm$.
\end{pf}

We finish this section with the rather well-known, important category-theoretic property of Banach spaces. Thanks to it, one can construct \Gurarii\ spaces.

\begin{lm}[The Pushout Lemma]
Let $Z, X, Y$ be Banach spaces,
let $\map iZX$ be an isometric embedding and let $\map fZY$ be an $\eps$-isometric embedding, where $\eps > 0$.
Then there exist a Banach space $W$, an isometric embedding $\map jYW$ and an $\eps$-isometric embedding $\map gXW$ for which the diagram
$$\xymatrix{
Y \ar[r]^j & W \\
Z \ar[u]^f \ar[r]_i & X \ar[u]_g
}$$
commutes.
Furthermore, if $X$, $Y$ are finite-dimensional then so is $W$.
\end{lm}

\begin{pf}
For simplicity, let us assume that $i$ is the inclusion $Z\subs X$.
Define $W = (X \oplus Y) \by \Delta$, where $X \oplus Y$ is endowed with the $\ell_1$ norm and
$$\Delta = \setof{ \pair z{-f(z)} }{z\in Z}.$$
Let $g$ and $j$ be the quotients of the canonical embeddings, i.e. $g(x) = \pair x0 + \Delta$ and $j(y) = \pair 0y + \Delta$ for $x\in X$, $y\in Y$.
Observe that
$$\norm {g(x)} = \inf_{z \in Z} \Bigl ( \norm{x + z}_X + \norm{ - f(z)}_Y \Bigr) \loe \norm{x}_X.$$
Similarly, 
$$\norm {j(y)} = \inf_{z\in Z} \Bigl( \norm{z}_X + \norm {y - f(z)}_Y \Bigr) \loe \norm{y}_Y.$$
It remains to estimate $\norm{g(x)}$ and $\norm{j(y)}$ from below.

Fix $x\in X$. Given $z\in Z$, we have
$$\norm {x + z}_X + \norm {-f(z)}_Y \goe (1+\eps)^{-1} \Bigl(\norm {x + z}_X + \norm {-z}_X\Bigr) \goe (1 + \eps)^{-1} \norm {x}_X.$$
It follows that $\norm{g(x)} \goe (1 + \eps)^{-1}\norm {x}_X$.

Now fix $y\in Y$. Given $z\in Z$, we have
$$\norm {z}_X  + \norm {y - f(z)}_Y \goe (1 + \eps)^{-1} \Bigl( \norm {f(z)}_Y + \norm {y - f(z)}_Y \Bigr) \goe \norm {y}_Y.$$
Thus $\norm{j(y)} \goe \norm {y}_Y$.
This completes the proof.
\end{pf}

We shall use several times the ``isometric" version of the Pushout Lemma:
Namely, if $f$ in the statement above is an isometric embedding then so is $g$.
Note also that the lemma above is valid when ``$\eps$-isometric embedding" is replaced by ``linear operator of norm $\loe 1+\eps$".
The proof is the same (see \cite{ACCGMud} for more details).

A word of explanation on the name ``Pushout Lemma" is in place. Namely, the commutative square from the lemma is usually called an \emph{amalgamation} of $X$ and $Y$ or, more precisely, of $i$ and $f$.
It turns out however that the amalgamation constructed in the proof is the pushout of $i$ and $f$ in the category of Banach spaces with linear operators of norm $\loe1$.
Specifically, given arbitrary bounded linear operators $\map TXV$, $\map SYV$ such that $T\cmp i = S \cmp f$, there exists a unique linear operator $\map hWV$ satisfying $h\cmp g = T$ and $h\cmp f = S$.
Finally, the norm of $h$ does not exceed $\max(\norm T, \norm S)$.

Recall that a space $Y \subs X$ is \emph{complemented} (more precisely: \emph{$k$-complemented}) in $X$ if there exists a projection $\map PXX$ of norm $\loe k$ and such that $Y = \im P$.
Officially, $P$ is a \emph{projection} if $P^2 = P$, however we shall say that a linear operator $\map PXY$ is a \emph{projection} if $Y\subs X$ and $P\rest Y = \id Y$.
It is clear that both definitions lead to the same concept.

Coming back to the previous remarks, the following property of a pushout deserves some attention:

\begin{lm}[cf. \cite{ACCGMud}]\label{Lpuszproj}
Under the assumptions of the Pushout Lemma, if $f$ is a linear operator of norm $\loe1$ and $\img iZ$ is $k$-complemented in $X$ then $\img jY$ is $k$-complemented in $W$.

Furthermore, if $i$ and $f$ are inclusions then every bounded projection from $X$ onto $Z$ extends to a projection from $W$ onto $Y$, preserving the norm.
\end{lm}

\begin{pf}
Let $\map PXZ$ be such that $P \cmp i = \id Z$ and $\norm P \loe k$.
Define $T = f\cmp P$ and $S = \id Y$.
Then $\norm T \loe k\norm f \loe k$, $T \cmp i = S \cmp f$, therefore by the property of the pushout, there exists a unique operator $\map hWZ$ of norm $\loe k$, such that $h\cmp g = T$ and $h\cmp j = S = \id Y$.
In particular, $j\cmp h$ gives a projection onto $\img jY \subs W$.
Finally, if $i$, $f$ are inclusions then $h\cmp g = T$ translates to $h \rest X = P$.
\end{pf}

Recall that a finite-dimensional Banach space $X$ is \emph{polyhedral} if its unit ball is a polyhedron.
In other words, there exist functionals $\phi_0,\phi_1,\dots,\phi_{k-1} \in X^*$ such that
$$\norm x = \max_{i<n} |\phi_i(x)|$$
for every $x\in X$.
An infinite-dimensional Banach space is called \emph{polyhedral} if each of its finite-dimensional subspaces is polyhedral.
Typical examples of polyhedral Banach spaces are $\ell_1(n)$ and $\ela(n)$, the $n$-dimensional variants of $\ell_1$ and $\ela$, respectively.
The spaces $\ela(n)$ play a special role, due to the following two facts.

\begin{prop}
A finite-dimensional Banach space is polyhedral if and only if it embeds isometrically into $\ela(n)$ for some $\Ntr$.
\end{prop}

\begin{pf}
Let the norm of $X$ be of the form
$$\norm x = \max_{i<n}|\phi_i(x)|,$$
where $\phi_0,\dots, \phi_{n-1} \in X^*$.
Define $\map eX{\ela(n)}$ by $e(x)(i) = \phi_i(x)$.
It is clear that $e$ is an isometric embedding.

Conversely, it is obvious that $\ela(n)$ is polyhedral, and every subspace of a polyhedral space is polyhedral.
\end{pf}

\begin{prop}
For every $\Ntr$ the space $\ela(n)$ is $1$-injective, that is, given a pair of Banach spaces $E \subs F$, every bounded linear operator $\map TE{\ela(n)}$ extends to a linear operator $\map {\til T}F{\ela(n)}$ so that $\norm {\til T} = \norm T$ holds.
\end{prop}

\begin{pf}
Fix $T$ and define $T_i(x) = T(x)(i)$ for $x\in E$.
By the Hahn-Banach Theorem, each $T_i$ extends to a linear functional $\til T_i$, preserving the norm.
Define $\til T(x)(i) = \til T_i(x)$ for $x\in F$. It is clear that $\til T$ extends $T$ and $\norm{\til T} = \norm T$.
\end{pf}

The proof of the following fact is an easy exercise, noticing that the norm of the pushout space is the convex hull of two polyhedra.

\begin{lm}\label{Lwielosciany}
Let $\map iZX$, $\map jZY$ be two isometric embeddings of finite-dimen\-sio\-nal polyhedral spaces. Then there exist a polyhedral space $W$ and isometric embeddings $\map {i'}XW$ and $\map {j'}YW$ such that the square
$$\xymatrix{
Y \ar[r]^{j'} & W \\
Z \ar[u]^j \ar[r]_i & X \ar[u]_{i'}
}$$
is commutative.
Furthermore, $W$ can be taken to be the space coming from the Pushout Lemma.
\end{lm}

\section{The separable \Gurarii\ space}\label{Segurademolinaaa}

This section has a survey character. We introduce the definition of a \Gurarii\ space, show its existence, uniqueness and basic properties.

\subsection{Universal disposition}

We start with some general definitions, originally due to \Gurarii~\cite{Gurarii}.

\begin{df}
Let $\fK$ be a class of Banach spaces (in most cases: either all finite-dimensional spaces or all separable spaces).
A Banach space $X$ is of \emph{(almost) universal disposition} for $\fK$ if for every pair of spaces $S\subs T$, both in $\fK$, for every isometric embedding $\map fSX$ (and for every $\eps > 0$), there exists an ($\eps$-)isometric embedding $\map gTX$ such that $g \rest S = f$.
If this holds, we shall write briefly ``$X$ is (almost) $\ud(\fK)$".
We shall write $\udfd$ and $\udsep$ for ``universal disposition for finite-dimensional spaces" and ``universal disposition for separable spaces", respectively.
\end{df}

\begin{df}
A Banach space is \emph{\Gurarii} if it is of almost universal disposition for finite-dimensional spaces.
A \emph{strong \Gurarii\ space} is a Banach space of universal disposition for finite-dimensional spaces.
\end{df}

The starting point of our study is the following result.

\begin{tw}[\Gurarii~\cite{Gurarii}]
There exists a separable \Gurarii\ space.
\end{tw}

We shall present two constructions in Subsection~\ref{SecPnfgvo}.

\subsection{Isometric uniqueness}

A standard back-and-forth argument shows that every two separable \Gurarii\ spaces are almost isometric. Below we sketch the arguments showing isometric uniqueness.

The following lemmas come from \cite{KubSol}. The proof of the first one is a bit technical, yet completely elementary. The second lemma follows directly from the first one, applying the definition of a \Gurarii\ space.

\begin{lm}\label{Lemkjperg}
Let $\map fXY$ be a strict $\eps$-isometric embedding of Banach spaces, $\eps > 0$.
Then there exist a Banach space $Z$ and isometric embeddings $\map gYZ$, $\map hXZ$, such that $\norm{g \cmp f - h} < \eps$.
\end{lm}

\begin{lm}\label{Llawofgurariii}
Let $G$ be a \Gurarii\ space. Then for every pair $X, Y$ of finite-dimensional Banach spaces such that $X \subs G$, for every $\eps > 0$, for every $\delta > 0$, for every strict $\eps$-isometric embedding $\map fXY$ there exists a $\delta$-isometric embedding $\map jYG$ such that $\norm{j f(x) - x} < \eps \norm x$ for every non-zero $x\in X$.
\end{lm}

\begin{tw}[Lusky~\cite{Lusky}]
Every two separable \Gurarii\ spaces are linearly isometric.
\end{tw}

\begin{pf}
Let $E$ and $F$ be two separable \Gurarii\ spaces.
Define inductively two sequences of linear operators $\map {f_n}{X_n}{Y_n}$ and $\map {g_n}{Y_n}{X_{n+1}}$ satisfying the following conditions.
\begin{enumerate}
	\item[(i)] $X_n\subs E$ and $Y_n\subs F$ are finite-dimensional spaces.
	\item[(ii)] $f_n$ and $g_n$ are $2^{-n}$-isometric embeddings.
	\item[(iii)] $\norm{g_n f_n(x) - x} < 2^{-n} \norm x$ for every $x\in X_n$.
	\item[(iv)] $\norm{f_{n+1} g_n(y) - y} < 2^{-n} \norm y$ for every $y\in Y_n$.
\end{enumerate}
We start with $X_0 = 0$ and we take $Y_0$ to be any finite-dimensional subspace of $F$.
We find $g_0$ by using Lemma~\ref{Llawofgurariii}.
Having defined $f_n$ and $g_n$, we use Lemma~\ref{Llawofgurariii} both for $E$ and $F$ to find first $f_{n+1}$ and next $g_{n+1}$.
Note that we have some freedom to choose the subspaces $X_{n+1}$ and $Y_{n+1}$.
Thus, the inductive construction can be carried out so that $\bigcup_{\ntr}X_n$ is dense in $E$ and $\bigcup_{\ntr}Y_n$ is dense in $F$.

Given $x\in X_n$, using (iv) and (ii), we have
$$\norm{f_n(x) - f_{n+1}g_n f_n(x)} < 2^{-n}\norm{f_n(x)} \loe 2^{-n+1}.$$
Similarly, using (ii) and (iii), we get
$$\norm{f_{n+1}(x) - f_{n+1}g_n f_n(x)} \loe \norm{f_{n+1}} \cdot \norm{x - g_n f_n(x)} < 2^{-n+1}.$$
Thus $\norm{f_n(x) - f_{n+1}(x)} < 2^{-n+2}$.
It follows that the sequence $\sett{f_n}{\ntr}$ is pointwise convergent.
Its limit extends uniquely to an isometry $\map {f_\infty}EF$.
The same arguments show that $\sett{g_n}{\ntr}$ pointwise converges to an isometry $\map {g_\infty}FE$.
Finally, (iii) and (iv) show that $g_\infty \cmp f_\infty = \id E$ and $f_\infty \cmp g_\infty = \id F$.
\end{pf}

From now on, we can speak about \emph{the} \Gurarii\ space, the unique separable space of almost universal disposition for finite-dimensional spaces.
This space will always be denoted by $\G$.

The proof above is actually a simplified version of that in~\cite{KubSol}, where it is shown that for every strict $\eps$-isometry $f$ between finite-dimensional subspaces of $\G$ there exists a bijective isometry $\map h\G\G$ such that $\norm{f-h} < \eps$.

\subsection{A criterion for being \Gurarii}

Note that there are continuum many isometric types of finite-dimensional Banach spaces.
Thus, to check that a given Banach space is \Gurarii, one needs to show the existence of suitable extensions of continuum many isometric embeddings.
Of course, this can be relaxed. One way to do it is to consider a natural countable subcategory of the category of all finite-dimensional Banach spaces.

We need to introduce some notation.
Every finite-dimensional Banach space $E$ is isometric to $\Err^n$ with some norm $\anorm$.
We shall say that $E$ is \emph{rational} if it is isometric to $\pair{\Err^n}{\anorm}$, such that the unit sphere is a polyhedron whose all vertices have rational coordinates.
Equivalently, $X$ is rational if, up to isometry, $X = \Err^n$ with a ``maximum norm" $\anorm$ induced by finitely many functionals $\phi_0,\dots, \phi_{m-1}$ such that $\img {\phi_i}{\Qyu^n} \subs \Qyu$ for every $i < m$. More precisely, $$\norm x = \max_{i < m}{|\phi_i(x)|}$$
for $x\in \Err^n$.
Typical examples of rational Banach spaces are $\ell_1(n)$ and $\ela(n)$, the $n$-dimensional variants of $\ell_1$ and $\ela$, respectively.
On the other hand, for $1 < p < \infty$, $n > 1$, the spaces $\ell_p(n)$ are not rational.
Of course, every rational Banach space is polyhedral.

It is clear that there are (up to isometry) only countably many rational Banach spaces and for every $\eps > 0$, every finite-dimensional space has an $\eps$-isometry onto some rational Banach space.

A pair of Banach spaces $\pair EF$ will be called \emph{rational} if $E \subs F$ and, up to isometry, $F = \Err^n$ with a rational norm, and $E\cap \Qyu^n$ is dense in $E$.
Note that if $\pair EF$ is a rational pair then both $E$ and $F$ are rational Banach spaces.
It is clear that there are, up to isometry, only countably many rational pairs of Banach spaces.

\begin{tw}\label{Twkrit}
Let $X$ be a Banach space. Then $X$ is \Gurarii\ if and only if it satisfies the following condition.
\begin{enumerate}
	\item[\rm{(G)}] Given $\eps > 0$, given a rational pair of spaces $\pair EF$, for every strict $\eps$-isometric embedding $\map fEX$ there exists an $\eps$-isometric embedding $\map gFX$ such that $$\norm {g\rest E - f} \loe \eps.$$
\end{enumerate}
Furthermore, in condition \rm{(G)} it suffices to consider $\eps$ from a given set $T\subs (0,+\infty)$ with $\inf T = 0$.
\end{tw}

\begin{pf}
Every \Gurarii\ space satisfies (G), almost by definition.
Assume $X$ satisfies (G).
Fix two finite-dimensional spaces $E\subs F$ and fix an isometric embedding $\map fEX$.
Fix $\eps > 0$.
Fix a linear basis $\Bee = \{ e_0, \dots, e_{m-1}\}$ in $F$ so that $\Bee \cap E = \{ e_0, \dots, e_{k-1}\}$ is a basis of $E$ (so $E$ is $k$-dimensional and $F$ is $m$-dimensional).
Choose $\delta > 0$ small enough.
In particular, $\delta$ should have the property that for every linear operators $\map{h,g}FX$, if
$\max_{i < m}\norm{h(e_i) - g(e_i)} < \delta$ then $\norm{h - g} < \eps/3$.
In fact, $\delta$ depends on the norm of $F$ only; a good estimation is $\eps/(3M)$, where
$$M = \sup\bigsetof{ \sum_{i < m}|\lam_i| }{ \bnorm{\sum_{i < m}\lam_i e_i} = 1}.$$
Now choose a $\delta$-equivalent norm $\anorm'$ on $F$ such that $E\subs F$ becomes a rational pair (in particular, the basis $\Bee$ gives a natural coordinate system in which all $e_i$s have rational coordinates).
The operator $f$ becomes a $\delta$-isometric embedding, therefore by (G) there exists a $\delta$-isometric embedding $\map gFX$ such that $\norm{f - g\rest E}' < \delta$.

Now let $\map hFX$ be the unique linear operator satisfying $h(e_i) = f(e_i)$ for $i < k$ and $h(e_i) = g(e_i)$ for $k \loe i < m$.
Then $h\rest \Bee$ is $\delta$-close to $g\rest \Bee$ with respect to the original norm, therefore $\norm{h - g} < \eps/3$.
Clearly, $h\rest E = f$.
If $\delta$ is small enough, we can be sure that $g$ is an $\eps/3$-isometric embedding with respect to the original norm of $F$.
Finally, assuming that $\eps < 1$, a standard calculation shows that $h$ is an $\eps$-isometric embedding, being $(\eps/3)$-close to $g$.

The ``furthermore" part obviously follows from the arguments above.
\end{pf}

Note that, for a given separable Banach space $X$, the criterion stated above can be applied by ``testing" countably many almost isometric embeddings, namely, only those that map rational vectors to a fixed countable dense subset of $X$.
More precisely, given a dense set $D\subs X$, every strict $\eps$-isometric embedding $\map f{\Err^n}X$ (where $\Err^n$ is endowed with some rational norm) can be approximated by strict $\eps$-isometric embeddings $\map g{\Err^n}X$ satisfying $\img g{\Qyu^n} \subs D$.

Theorem~\ref{Twkrit} together with Lemma~\ref{Lemkjperg} provide another natural criterion for being \Gurarii.

\begin{tw}\label{TcirtwjrF}
A Banach space $X$ is \Gurarii\ if and only if it satisfies the following condition.
\begin{enumerate}
	\item[\rm{(F)}] Given $\eps, \delta > 0$, given a rational pair of spaces $\pair EF$, for every strict $\eps$-isometric embedding $\map fEX$ there exists a $\delta$-isometric embedding $\map gFX$ such that $\norm{f - g\rest E} < \eps$.
\end{enumerate}
\end{tw}

\begin{pf}
It is clear that (F) implies (G) and, by Theorem~\ref{Twkrit} this implies that $X$ is \Gurarii.
It remains to show that every \Gurarii\ space satisfies (F).
For this aim, fix a rational pair $\pair EF$ and a strict $\eps$-isometric embedding $\map fEX$.
Let $Y = \img fE$.
By Lemma~\ref{Lemkjperg}, there are a finite-dimensional space $Z$ and isometric embeddings $\map iEZ$ and $\map jYZ$ such that $\norm{ j\cmp f - i} < \eps$.
Using the Pushout Lemma, we can extend $Z$ so that it also contains $F$.
Since $X$ is \Gurarii, there exists a $\delta$-isometric embedding $\map hZX$ extending $j^{-1}$.
Finally, $g = h \rest F$ is as required.
\end{pf}

\subsection{Two constructions}\label{SecPnfgvo}

There are several ways to see the existence of the \Gurarii\ space $\G$.
Actually, in Theorem~\ref{Tgnoir} below, we shall show the existence of strong \Gurarii\ spaces; in view of Theorem~\ref{Tgcerugh} below, such spaces contain many isometric copies of the \Gurarii\ space.
However, this is a rather indirect way of showing the existence of $\G$.
A direct way is to construct a certain chain of finite-dimensional spaces.
The crucial point is the Pushout Lemma.

\begin{tw}[\Gurarii\ \cite{Gurarii}, Gevorkjan \cite{Gevork}]
The \Gurarii\ space exists and is isometrically universal for all separable Banach spaces. 
\end{tw}

\begin{pf}
Fix a separable Banach space $X$ and fix a countable dense set $D\subs X$.
Fix a rational pair of Banach spaces $E\subs F$, fix a linear basis $B$ in $E$ consisting of vectors with rational coordinates, and fix $\eps > 0$.
Furthermore, fix a strict $\eps$-isometric embedding $\map fEX$ such that $\img fB \subs D$.
Using the Pushout Lemma, we can find a separable Banach space $X' \sups X$ such that $f$ extends to a strict $\eps$-isometric embedding $\map gF{X'}$.
Note that there are only countably many pairs of rational Banach spaces and almost isometric embeddings as described above.
Thus, there exists a separable Banach space $G(X) \sups X$ such that, given a rational pair $E\subs F$, for every $\eps$-isometric embedding $\map fEX$ there exists an $\eps$-isometric embedding $\map gFX$ such that $g\rest E$ is arbitrarily close to $f$.

Repeat this construction infinitely many times.
Namely, let $G = \cl( \bigcup_{\ntr} X_n )$, where $X_0 = X$ and $X_{n+1} = G(X_n)$ for $\ntr$.
Clearly, $G$ is a separable Banach space.
By Theorem~\ref{Twkrit}, $G$ is the \Gurarii\ space.

Since the space $X$ was chosen arbitrarily, this also shows that the \Gurarii\ space contains an isometric copy of every separable Banach space.
\end{pf}

Next we show how to construct the \Gurarii\ space as a ``random" or ``generic" Banach space. Uncountable variants, forcing the universe of set theory to be extended, have been recently studied by Lopez-Abad and Todorcevic~\cite{JordiS}.
Our idea is similar in spirit to that of \Gurarii\ from \cite{Gurarii}, however it does not use any topological structure on spaces of norms.

Recall that $c_{00}$ denotes the linear subspace of $\Err^\omega$ consisting of all vectors with finite support.
In other words, $x\in c_{00}$ iff $x\in \Err^\nat$ and $x(n)=0$ for all but finitely many $\ntr$.
Given a finite set $S\subs\nat$, we shall identify each space $\Err^n$ with the suitable subset of $c_{00}$.

Let $\poset$ be the following partially ordered set.
An element of $\poset$ is a pair $p = \pair {S_p}{\anorm_p}$, where $S_p\subs\nat$ is a finite set and ${\anorm_p}$ is a norm on $\Err^{S_p} \subs c_{00}$.
We define $p \loe q$ iff $S_p \subs S_q$ and $\anorm_q$ extends $\anorm_p$.
Clearly, $\poset$ is a partially ordered set.
Suppose $$p_0 < p_1 < p_2 < \dots$$
is a sequence in $\poset$ such that the chain of sets $\bigcup_{\ntr}S_{p_n} = \nat$.
Then $c_{00}$ naturally becomes a normed space.
Let $X$ be the completion of $c_{00}$ endowed with this norm.
We shall call it the \emph{limit} of $\ciag p$ and write $X = \lim_{n\to\infty} p_n$.
It is rather clear that every separable Banach space is of the form $\lim_{n\to\infty}p_n$ for some sequence $\ciag p$ in $\poset$.
We are going to show that for a ``typical" sequence in $\poset$, its limit is the \Gurarii\ space.

Given a partially ordered set $\poset$, recall that a subset $D\subs \poset$ is \emph{cofinal} if for every $p\in \poset$ there exists $d\in D$ with $p \loe d$.
Below is a variant of the well-known Rasiowa-Sikorski Lemma, which is actually an abstract version of the Baire Category Theorem.

\begin{lm}\label{RasSik}
Let $\poset$ be a partially ordered set and let $\Dee$ be a countable family of cofinal subsets of $\poset$. Then there exists a sequence
$$p_0 \loe p_1 \loe p_2 \loe \dots$$
such that for each $D\in\Dee$ the set $\setof{\ntr}{p_n \in D}$ is infinite.
\end{lm}

\begin{pf}
Let $\Dee = \setof{D_n}{\ntr}$ so that for each $D\in\Dee$ the set $\setof{\ntr}{D_n = D}$ is infinite. Using the fact that each $D_n$ is cofinal, construct inductively $\ciag p$ so that $p_n \in D_n$ for $\ntr$.
\end{pf}

A sequence $\ciag p$ satisfying the assertion of the lemma above is often called \emph{$\Dee$-generic}.

We now define a countable family of open cofinal sets which is good enough for producing the \Gurarii\ space.
Namely, fix a rational pair of spaces $\pair EF$, fix a positive integer $n$ and fix a rational embedding $\map fE{c_{00}}$, that is, an injective linear operator mapping vectors with rational coordinates to $c_{00}\cap \Qyu^{\nat}$.
The point is that there are only countably many possibilities for $E,f$.
Define
$D_{E,F,f,n}$
to be the set of all $p\in \poset$ such that $n\in S_p$ and $p$ satisfies the following implication:
If $f$ is a $(1/n)$-isometric embedding into $\pair {\Err^{S_p}}{\anorm_p}$, then there exists a $(1/n)$-isometric embedding $\map gF{\pair {\Err^{S_p}}{\anorm_p}}$ such that $g \rest E = f$.

\begin{claim}
The set $D_{EF,f,n}$ is cofinal in $\poset$.
\end{claim}

\begin{pf}
Fix $p\in \poset$.
Suppose that $f$ is a $(1/n)$-isometric embedding into $\pair{\Err^{S_p}}{\anorm_p}$ (otherwise clearly $p\in D_{E,F,f,n}$).
Using the Pushout Lemma, find a finite-dimensional Banach space $W$ extending $\pair{\Err^{S_p}}{\anorm_p}$ and a $(1/n)$-isometric embedding $\map gFW$ such that $g \rest F = f$.
We may assume that $W = \pair{\Err^T}{\anorm_W}$ for some $T\sups S_p$, where the norm $\anorm_W$ extends $\anorm_p$.
Let $q = \pair T{\anorm_W} \in \poset$.
Clearly, $p \loe q$ and $q\in D_{E,F,f,n}$.
\end{pf}

Let $\Dee$ consist of all sets of the form $D_{E,F,f,n}$ as above.
Then $\Dee$ is countable, therefore applying Lemma~\ref{RasSik} we obtain a $\Dee$-generic sequence $\ciag p$.

\begin{tw}\label{Tzbirczporz}
Let $\Dee$ be as above and let $\ciag p$ be a $\Dee$-generic sequence.
Then the space $\lim_{n\to\infty}p_n$ is \Gurarii.
\end{tw}

\begin{pf}
Let $X = \lim_{n\to\infty}p_n$.
Notice that $\bigcup_{\ntr}S_{p_n}=\nat$.
Fix a positive integer $k$, fix a rational pair of spaces $\pair EF$ and fix a $1/(k+1)$-isometric embedding $\map fEX$.
We can modify $f$ in such a way that it remains to be a $(1/k)$-isometric embedding, while at the same time $\img fE \subs c_{00}$, and it maps rational vectors into $c_{00}\cap \Qyu^\nat$.
Now $D_{E,F,f,k} \in \Dee$ therefore there exists $\ntr$ such that $p_n \in D_{E,F,f,k}$ and $\Err^{S_{p_n}}$ contains the range of $f$.
By the definition of $D_{E,F,f,k}$, $f$ extends to a $(1/k)$-isometric embedding $\map gFX$.
By Theorem~\ref{Twkrit}, this shows that $X$ is \Gurarii.
\end{pf}

Let us remark that some modifications of the poset $\poset$ still give the \Gurarii\ space. For instance, we can consider only polyhedral norms for $\anorm_p$, because the Pushout Lemma holds for this class.
We shall use this observation later.

\subsection{Schauder bases and Lindenstrauss spaces}

We now present the proof that the \Gurarii\ space has a monotone Schauder basis.
This fact has already been noticed by \Gurarii\ in \cite{Gurarii}.

Recall that a \emph{Schauder basis} in a separable Banach space $X$ is a sequence $\ciag e$ of non-zero vectors of $X$, such that for every $x\in X$ there exist uniquely determined scalars $\ciag \lam$ satisfying
$$x = \sum_{\Ntr}\lam_n e_n.$$
The series above is supposed to converge in the norm.
Given a Schauder basis $\ciag e$, one always has the associated \emph{canonical projections}
$$P_N \Bigl(\sum_{\Ntr}\lam_n e_n\Bigr) = \sum_{n < N}\lam_n e_n.$$
Note that each $P_N$ is a projection and $P_N P_M = P_{\min(N,M)}$ for every $N,M \in \Nat$.
By the Banach-Steinhaus principle, $\sup_{N \in\Nat}\norm {P_N} < +\infty$.
The basis is \emph{monotone} if $\norm{P_N}\loe1$ for each $N\in\Nat$.
We shall consider monotone Schauder bases only.
It turns out that the existence of a monotone Schauder basis can be deduced from the canonical projections:

\begin{prop}[Mazur]\label{Pbasia}
Let $X$ be a Banach space and let $\ciag P$ be a sequence of norm one projections such that $P_0 = 0$, $\dim(P_{n+1}X \by P_nX) \loe1$ and $P_n P_m = P_{\min(n,m)}$ for every $n,m \in \Nat$.
Then there exists a monotone Schauder basis $\ciag e$ in $X$ such that $\ciag P$ is  the sequence of canonical projections associated to $\ciag e$.
\end{prop}

\begin{pf}
Let us first prove that $\lim_{n\to\infty}P_n x = x$ for every $x\in X$.
For this aim, fix $x\in S_X$ and $\eps > 0$.
Find $n_0$ such that $\norm {x - y} < \eps/2$ for some $y\in P_{n_0}X$.
Given $n \goe n_0$, we have
$$\norm {P_n x - x} = \norm{P_n x - y} + \norm {y  - x} < \norm{P_n (x - y)} + \eps/2 < \eps.$$
Now let $\phi_n$ be such that $P_{n+1} x - P_n x = \phi_n(x) e_n$ for some $e_n \in S_X$.
Here we use the fact that $P_{n+1}X = P_n X \oplus \Err e_n$ for some $e_n\in \ker {P_n} \cap S_X$.
Finally, given $x\in X$, we have
$$x = \lim_{n\to\infty}P_n x = \lim_{N\to\infty}\sum_{n<N}(P_{n+1} - P_n)x = \lim_{N\to\infty}\sum_{n<N}\phi_n(x) e_n.$$
Finally, if $0 = \sum_{\Ntr}\lam_n e_n$ then, by easy induction, we show that $\lam_n = 0$ for every $\Ntr$. This shows that $\ciag e$ is a Schauder basis.
Clearly, $P_n$s are the canonical projections, therefore the basis is monotone.
\end{pf}

There is an important class of Banach spaces with monotone Schauder bases, containing the \Gurarii\ space:

\begin{df}
A Banach space $X$ is called a \emph{Lindenstrauss space} if it contains a directed (by inclusion) family $\Ef$ such that $\bigcup\Ef$ is dense in $X$ and each $F\in\Ef$ is linearly isometric to some $\ela(n)$.
\end{df}

Lindenstrauss spaces were introduced under the name \emph{$(\pi_1^\infty)$ spaces} by Michael \& Pe{\l}\-czy\'n\-ski~\cite{MicPel}.
It easily follows from their results that a space $X$ is Lindenstrauss if and only if it is \emph{almost 1-injective} for finite-dimensional spaces,
that is, given a pair of finite-dimensional spaces $R\subs S$ and a norm one operator $\map fRX$, there exists for every $\eps>0$ an extension $\map gSX$ of $f$ satisfying $\norm g \loe 1+\eps$.
Recall that every space of the form $\ela(n)$ is $1$-injective for all Banach spaces.
It is worth pointing out that every Lindenstrauss space contains a dense linear subspace that is $1$-injective for finite-dimensional spaces.
The basic infinite-dimensional example of a Lindenstrauss space is $c_0$. This particular space is $1$-injective for finite-dimensional spaces.

\begin{tw}[Michael \& Pe{\l}czy\'nski \cite{MicPel}]\label{Tmicpelghriuw}
A separable infinite-dimensional Banach space is a Lindenstrauss space if and only if it is the completion of the union of a chain $\bigcup_{\Ntr}X_n$, where each $X_n$ is linearly isometric to $\ela(n)$.
\end{tw}

This is an immediate consequence of the following interesting geometric property of finite-dimensional $\ela$ spaces:

\begin{lm}[Michael \& Pe{\l}czy\'nski \cite{MicPel}]\label{LMPgeometiras}
Given an isometric embedding $\map f{\ela(k)}{\ela(l)}$ with $k + 1 < l$, there exists a space $Z$ linearly isometric to $\ela(k+1)$ and such that $\img f{\ela(k)} \subs Z \subs \ela(l)$.
\end{lm}

Theorem~\ref{Tmicpelghriuw} combined with Proposition~\ref{Pbasia} gives the following

\begin{wn}[\Gurarii~\cite{Gur2}, Michael \& Pe{\l}czy\'nski~\cite{MicPel}]
Every separable Lindenstrauss space has a monotone Schauder basis.
\end{wn}

\begin{tw}[\Gurarii\ \cite{Gurarii}]\label{TGuryiLasy}
The \Gurarii\ space is a Lindenstrauss space.
\end{tw}

\begin{pf}
Using the Pushout Lemma, it is straightforward to see that $\G$ is almost 1-injective for finite-dimensional spaces, therefore it is Lindenstrauss.
Below we present a direct argument.

Let $\poset$ be the partially ordered set defined before Theorem~\ref{Tzbirczporz}.
Define $\poset_0$ to be the set of all $p\in\poset$ such that the norm $\anorm_p$ is polyhedral.
It is easy to verify that, with the same family $\Dee$ of cofinal sets, the limit of a $\Dee$-generic sequence is \Gurarii.
In fact, the only difference is in using the polyhedral variant of the Pushout Lemma, namely, Lemma~\ref{Lwielosciany}.
Now add to the family $\Dee$ the following set:
$$E = \setof{p\in\poset_0}{\pair {\Err^{S_p}}{\anorm_p} \text{ is linearly isometric to some }\ela (n)}.$$
Since all the norms $\anorm_p$ are polyhedral, the set $E$ is cofinal in $\poset_0$.
The limit of a $(\Dee\cup\sn E)$-generic sequence is necessarily a Lindenstrauss space; since such a sequence is also $\Dee$-generic, its limit is the \Gurarii\ space.
\end{pf}

\begin{tw}[Wojtaszczyk \cite{Wojt}]\label{TWojtaszz}
A separable Banach space is a Lindenstrauss space if and only if it embeds isometrically onto a 1-complemented subspace of $\G$.
\end{tw}

\begin{pf}
Note that, by the Pushout Lemma, the \Gurarii\ space is almost 1-injective for finite-dimensional spaces, therefore so is any of its 1-complemented subspaces.
This shows the ``only if" part of the theorem above.

Fix a separable Lindenstrauss space $X$ and let $\ciag X$ be a chain of spaces such that $X_0 = \sn0$, $\bigcup_{\Ntr}X_n$ is dense in $X$, and each $X_n$ is linearly isometric to some $\ela(k_n)$. By Theorem~\ref{Tmicpelghriuw}, we may assume that $k_n=n$ unless $X$ is finite-dimensional, where we simply put $X_n = X = \ela(\dim X)$.

Let us look back at the simple proof of Lemma~\ref{RasSik}, where a $\Dee$-generic sequence was constructed in the poset $\poset$ defined just before Theorem~\ref{Tzbirczporz} and $\Dee$ is the same countable collection of cofinal sets.
For convenience, we shall write $U(p)$ for the Banach space $\pair{\Err^{S_p}}{\anorm_p}$, where $p\in\poset$.

We claim that there exists a $\Dee$-cofinal sequence $\ciag p$ together with isometric embeddings $\map {i_n}{X_n}{U(p_n)}$ and norm one operators $\map {P_n}{U(p_n)}{X_n}$ such that $P_n \cmp i_n = \id{X_n}$, $i_{n+1}$ extends $i_n$ and $P_{n+1}$ extends $P_n$ for each $\Ntr$.
Recall that $\Dee$ was enumerated as $\ciag D$, so that each $D\in\Dee$ occurs infinitely many times.
Suppose $p_n$, $i_n$ and $P_n$ have been defined.
Using the Pushout Lemma, find $q \goe p_n$ and an isometric embedding $\map j{X_{n+1}}{U(q)}$ extending $i_n$.
The property of the pushout gives a norm one projection $\map Q{U(q)}{X_{n+1}}$ extending $P_n$ (see Lemma~\ref{Lpuszproj}).

Now, using the fact that $D_{n+1}$ is cofinal, find $p_{n+1}\in D_{n+1}$ so that $p_{n+1} \goe q$.
Finally, $i_{n+1} = j$, treated as an embedding into $U(p_{n+1})$ and $P_{n+1}$ is any extension of $Q$ preserving the norm, which exists because $X_{n+1}$ is linearly isometric to some $\ela(m)$.

This finishes the inductive construction.
By Theorem~\ref{Tzbirczporz}, we know that $\lim_{n\to\infty}p_n = \G$ and taking the pointwise limits of $i_n$ and $P_n$
we obtain an isometric embedding $\map iX\G$ and a norm one operator $\map P{\G}X$ such that $P \cmp i = \id X$. This shows that $\img iX$ is 1-complemented in $\G$.
\end{pf}

The arguments from the first part of the proof above can actually be viewed as another proof of the isometric universality of the \Gurarii\ space, since the embeddings $i_n$ can be constructed without assuming that $X$ is a Lindenstrauss space.

\section{Non-separable \Gurarii\ spaces}\label{Secnonsepeirge}

In this section we give a characterization of \Gurarii\ spaces in terms of skeletons.

Let $X$ be a Banach space.
A family $\Ef$ of closed linear subspaces of $X$ will be called a \emph{skeleton} in $X$ if the following conditions are satisfied.
\begin{enumerate}
	\item[(1)] Each $F\in \Ef$ is separable.
	\item[(2)] $\bigcup\Ef = X$.
	\item[(3)] For every $F_0,F_1\in \Ef$ there exists $G \in \Ef$ such that $F_0\cup F_1 \subs G$.
	\item[(4)] $\cl(\bigcup_{\ntr}F_n) \in \Ef$, whenever $\sett{F_n}{\ntr}$ is a countable chain in $\Ef$.
\end{enumerate}
A family of sets $\Ef$ satisfying (3) is called \emph{directed}.
The notion of a skeleton makes sense in non-separable Banach spaces, since $\Ef = \sn X$ is a skeleton if $X$ is separable.
Actually, notice that if $\Ef$ is a skeleton in $X$ then for every separable subset $S\subs X$ there exists $F\in\Ef$ satisfying $S \subs F$.
In particular, every skeleton in a separable Banach space $X$ must contain $X$.
The value of skeletons lies in the following well-known property.

\begin{prop}\label{Pskelegeres}
Let $\Ef$ and $\Gee$ be skeletons in a fixed Banach space $X$. Then $\Ef \cap \Gee$ is also a skeleton in $X$.
\end{prop}

\begin{pf}
It is clear that $\Ef \cap \Gee$ satisfies (1) and (4).
In order to prove (2) and (3) it suffices to show that for every separable subspace $S\subs X$ there exists $H\in \Ef \cap \Gee$ such that $S \subs H$.

Fix a separable set $S\subs X$. By the remark above, there exists $F_0 \in \Ef$ such that $S\subs F_0$.
Similarly, there exists $G_0 \in \Gee$ such that $F_0\subs G_0$.
By induction, we construct two increasing sequences $\ciag F$ and $\ciag G$ in $\Ef$ and $\Gee$ respectively, so that $F_n \subs G_n \subs F_{n+1}$ holds for every $\ntr$.
Finally, notice that $H = \cl(\bigcup_{\ntr}F_n) = \cl(\bigcup_{\ntr}G_n)$ belongs to both $\Ef$ and $\Gee$.
\end{pf}

We now turn to the announced characterization of \Gurarii\ spaces in terms of skeletons.

\begin{lm}\label{Lifhtgo}
Let $X$ be a \Gurarii\ space and let $S\subs X$ be a countable set.
Then there exists a subspace $Y\subs X$ linearly isometric to $\G$ and such that $S \subs Y$.
\end{lm}

\begin{pf}
This is a standard closing-off argument.
The criterion for being \Gurarii\ (Theorem~\ref{Twkrit}) actually requires checking countably many almost isometric embeddings.
The first step is to show that given a separable subspace $Z\subs X$ there exists a separable space $E(Z)$ (not uniquely determined) such that $Z \subs E(Z) \subs X$ and the following condition is satisfied.
\begin{enumerate}
\item[($\dagger$)] For every rational pair of spaces $\pair EF$, for every $\eps > 0$, for every strict $\eps$-isometric embedding $\map fEZ$ there exists a strict $\eps$-isometric embedding $\map gF{E(Z)}$ such that $\norm{f - g\rest E} < \eps$.
\end{enumerate}
Once we have proved this, we construct a chain of separable spaces $Z_0 \subs Z_1 \subs \dots \subs X$ such that $S\subs Z_0$, and $Z_{n+1} = E(Z_n)$ for every $\Ntr$.
Then, using Theorem~\ref{Twkrit}, we conclude that the space $Y = \cl(\bigcup_{\ntr}Z_n)$ is \Gurarii, because of condition ($\dagger$).
It remains to show the existence of $E(Z)$ satisfying ($\dagger$).

Fix $Z$ and fix a countable dense subset $D$ of $Z$.
Let $\Aaa$ consist of all quadruples of the form $\seq{E,F,f,\eps}$, where $E\subs F$ is a rational pair of (finite-dimensional) spaces, $\eps > 0$ is a rational number, and $\map fEZ$ is a strict $\eps$-isometric embedding such that $\img fB \subs D$, where $B$ is a fixed linear basis of $E$ consisting of vectors with rational coordinates.
These assumptions ensure us that $\Aaa$ is countable.

Using the fact that $X$ is \Gurarii, given $q = \seq{E,F,f,\eps}$, we know that there exists a strict $\eps$-isometric embedding $\map gFX$ such that $\norm{f - g\rest E} < \eps$.
Denote by $R_q$ the range of $g$.
Finally, take $E(Z)$ to be the closure of the union $Z \cup \bigcup_{q\in\Aaa}R_q$.
It is clear that ($\dagger$) is satisfied.
\end{pf}

\begin{lm}\label{Lkuqqrs}
Let $\ciag X$ be a chain of subspaces of a Banach space $X$ such that $X = \cl(\bigcup_{\ntr}X_n)$ and each $X_n$ is linearly isometric to $\G$.
Then $X$ is linearly isometric to the \Gurarii\ space $\G$.
\end{lm}

\begin{pf}
Fix finite-dimensional spaces $E \subs F$ and an isometric embedding $\map fEX$.
Fix $\eps > 0$.
Choose a linear map $\map gEX$ that is $\eps$-close to $f$ so that it is a strict $\eps$-isometric embedding and $\img gE \subs X_n$ for some $\ntr$.
Now, using the property of the \Gurarii\ space $X_n$, there exists an extension $\map hF{X_n}$ of $g$, that is also a strict $\eps$-isometric embedding.
Finally, $h\rest E$ is $\eps$-close to $f$. By Theorem~\ref{Twkrit}, this shows that $X$ is \Gurarii.
\end{pf}

\begin{tw}\label{Tgcerugh}
Let $X$ be a Banach space. The following properties are equivalent.
\begin{enumerate}
	\item[(a)] $X$ is a \Gurarii\ space.
	\item[(b)] $X$ has a skeleton consisting of subspaces isometric to the \Gurarii\ space $\G$.
	\item[(c)] There exists a directed family $\Gee$ of spaces isometric to $\G$, such that $\bigcup\Gee = X$.
\end{enumerate}
\end{tw}

\begin{pf}
(a)$\implies$(c)
Let $\Gee$ be the family of all subspaces of $X$ that are isometric to $\G$.
By Lemma~\ref{Lifhtgo}, $\bigcup\Gee = X$ and $\Gee$ is directed.
In fact, this follows from a stronger property of $\Gee$: every separable subset is covered by an element of $\Gee$.

(c)$\implies$(b)
Let $\Gee$ be as in (c) and let $\Ef$ be the family of all subspaces of $X$ that are isometric to $\G$.
We claim that $\Ef$ is a skeleton.
Condition (1) is obvious, (2) follows from the property of $\Gee$ and (4) follows from Lemma~\ref{Lkuqqrs}.
In order to show (3), it suffices to prove that every countable subset of $X$ is covered by an element of $\Ef$.
Fix $D = \setof{d_n}{\ntr} \subs X$ and, using directedness, construct inductively $G_0 \subs G_1 \subs \dots$ in $\Gee$ so that $d_n \in G_n$.
Then $F = \cl(\bigcup_{\ntr}G_n)$ is an element of $\Ef$ and $D \subs F$.

(b)$\implies$(a)
Fix two finite-dimensional spaces $A \subs B$ and an isometric embedding $\map fAX$. Then $\img fA$ is finite-dimensional, therefore there exists $F\in \Ef$ such that $\img fA \subs F$. Since $F$ is the \Gurarii\ space, given any $\eps>0$, $f$ can be extended to an $\eps$-isometry $\map gBF$.
\end{pf}

The following corollary improves \cite[Thm. 6.1]{ACCGMud}, where the same was shown for Banach spaces of universal disposition for separable spaces.

\begin{wn}\label{wnhorig}
No complemented subspace of a $C(K)$ space (or, more generally, an M-space) can be \Gurarii.
\end{wn}

\begin{pf}
Suppose $X\subs C(K)$ is a \Gurarii\ space and $\map P{C(K)}X$ is a projection.
Let $\Ef$ be a skeleton in $C(K)$ consisting of spaces of continuous functions over some metric compacta.
By Theorem~\ref{Tgcerugh}, $X$ has a skeleton $\Gee$ such that each $G\in\Gee$ is isometric to the \Gurarii\ space $\G$.
A standard closing-off argument (see the proof of~\cite[Thm. 6.1]{ACCGMud}) shows that there are $F\in\Ef$ and $G\in \Gee$ such that $PF = G$.
The final contradiction comes from \cite[Cor. 5.4]{ACCGMud}, saying that the \Gurarii\ space is not complemented in any C(K) space.

The arguments above can be repeated when C(K) spaces are replaced by M-spaces (see the comments in Sections 5,6 of \cite{ACCGMud}).
\end{pf}

It should be noted that Corollary~\ref{wnhorig} can actually be derived from \cite[Thm. 6.1]{ACCGMud}, using another result from \cite{ACCGMud} saying that ultraproducts of \Gurarii\ spaces are $\udsep$, while ultraproducts of $C(K)$ spaces are again $C(K)$ spaces.
However, our argument using skeletons is elementary and perhaps more illustrative.

\begin{tw}\label{Tisaisaes}
Every Banach space embeds isometrically into a \Gurarii\ space of the same density.
\end{tw}

\begin{pf}
We use induction on the density of the space. The statement is true for separable spaces, so fix a cardinal $\kappa > \aleph_0$ and suppose the statement holds for Banach spaces of density $< \kappa$.

Fix a Banach space $X$ of density $\kappa$.
Then $X$ is the completion of the union of a chain $\sett{X_\al}{\al<\kappa}$ starting from a separable space $X_0$ and such that $\dens {X_\al} < \kappa$ for every $\al < \kappa$.
We may assume that this chain is continuous, i.e., $X_\delta$ is the closure of $\bigcup_{\xi < \delta}X_\xi$, whenever $\delta$ is a limit ordinal.
We construct a sequence of isometric embeddings $\map {f_\al}{X_\al}{G_\al}$, where each $G_\al$ is a \Gurarii\ space of density $<\kappa$, and $G_\al \subs G_\beta$, $f_\beta \rest G_\al = f_\al$ whenever $\al < \beta$.

Suppose $G_\al$ and $f_\al$ have been constructed for $\al < \eta$.
If $\eta$ is a limit ordinal, we take $G_\eta$ to be the completion of $\bigcup_{\xi < \eta}G_\xi$.
By Theorem~\ref{Tgcerugh}, we know that $G_\eta$ is a \Gurarii\ space.
The embedding $f_\eta$ is uniquely determined.

Now suppose $\eta = \beta + 1$. Using the Pushout Lemma, we find a space $W \sups G_\beta$ so that $f_\beta$ extends to an isometric embedding $\map j{X_{\beta+1}}W$.
Note that $\dens W < \kappa$.
Using the inductive hypothesis, there exists a \Gurarii\ space $G_{\beta+1} \sups W$ such that $\dens {G_{\beta+1}} = \dens W$. We define $f_{\beta+1} = j$.

Finally, the sequence $\sett{f_\al}{\al<\kappa}$ determines an isometric embedding of $X$ into $G = \cl(\bigcup_{\al<\kappa}G_\al)$.
Clearly, $\dens G = \kappa$ and $G$ is \Gurarii\ by Theorem~\ref{Tgcerugh}.
\end{pf}

It seems that there are many non-isomorphic \Gurarii\ spaces of density $\aleph_1$.
We show that some of them have many projections.
Recall that a \emph{projectional resolution of the identity} (briefly: \emph{PRI}) in a Banach space is a transfinite sequence of norm one projections $\sett{P_\al}{\al<\omega_1}$ whose images are separable, form a continuous chain covering the space, and $P_\al P_\beta = P_{\min\{\al,\beta\}}$ holds for every $\al,\beta < \omega_1$.
The notion of a PRI is usually defined for arbitrary non-separable Banach spaces, see \cite{DGZ} and \cite {Fab} for more information.
It seems that PRI is the main tool for proving certain properties of a non-separable Banach space by transfinite induction.
For example, every Banach space of density $\aleph_1$ with a PRI admits a bounded one-to-one linear operator into $c_0(\omega_1)$ (see, e.g., \cite[Cor. 17.5]{KKLP}).

\begin{tw}\label{Tprispris}
There exists a \Gurarii\ space $E$ of density $\aleph_1$ that has a projectional resolution of the identity.
\end{tw}

\begin{pf}
First of all, there exists a norm one projection $\map Q\G\G$ such that $\ker Q$ is non-trivial.
This follows immediately from the proof of Theorem~\ref{TWojtaszz}, where we can at the first step ensure that the embedding of $\G$ into $\G$ is not the identity.

We now construct a continuous chain of separable spaces $\sett{G_\al}{\al < \omega_1}$ with the following properties.
\begin{enumerate}
	\item[(i)] Each $G_\al$ is linearly isometric to $\G$.
	\item[(ii)] For each $\al < \omega_1$ there exists a projection $\map {Q^{\al+1}_\al}{G_{\al+1}}{G_\al}$, isometric to $Q$.
\end{enumerate}
Property (ii) ensures us that the chain is strictly increasing and its union $G_{\omega_1} = \bigcup_{\al < \omega_1}G_\al$ is indeed of density $\aleph_1$.

By Theorem~\ref{Tgcerugh}, $G_{\omega_1}$ is a \Gurarii\ space and by~\cite{Kub_lines} (see also~\cite[Thm. 17.5]{KKLP}) it has a projectional resolution of the identity.
\end{pf}

Note that there are Banach spaces of density $\aleph_1$, not embeddable into any Banach space with a PRI (e.g., spaces with uncomplemented copies of $c_0$, see Section~\ref{Sroleceezero} below).
Thus, by Theorems~\ref{Tisaisaes} and \ref{Tprispris}, there are at least two non-isomorphic \Gurarii\ spaces of density $\aleph_1$.

\section{Spaces of universal disposition for larger classes}\label{SecUDudud}

In this section we discuss spaces of $\ud(\fD_{<\kappa})$, where $\fD_{<\kappa}$ is the class of Banach spaces of density $< \kappa$.

Recall that a Banach space is \emph{isometrically universal} for a class $\fK$ of spaces, if it contains an isometric copy of every space from $\fK$.
The following general fact is well-known, we state it for the sake of completeness.
A special case (for $\kappa = \aleph_0$) is contained in \cite{Gurarii}.

\begin{prop}\label{Prbeou}
Let $\kappa$ be an infinite regular cardinal and denote by $\fD_{<\kappa}$ the class of all Banach spaces of density $< \kappa$. If $\kappa = \aleph_0$, let $\fD_{<\kappa}$ be the class of all finite-dimensional Banach spaces.
\begin{enumerate}
	\item[(0)] Let $U$ be a Banach space of $\ud(\fD_{<\kappa})$.
Then for every pair of spaces $X\subs Y$ such that $\dens X < \kappa$ and $\dens Y \loe \kappa$, every isometric embedding $\map fXU$ extends to an isometric embedding $\map gYU$.
	\item[(1)] Every Banach space of $\ud(\fD_{<\kappa})$ is isometrically universal for the class of Banach spaces of density $\loe \kappa$.
	\item[(2)] Assume $\kappa$ is a regular cardinal and let $U$, $V$ be two Banach spaces of $\ud(\fD_{<\kappa})$ and of density $\kappa$. Then every linear isometry $\map fXY$ such that $X\subs U$, $Y\subs V$ and $X,Y \in \fD_{<\kappa}$, extends to a bijective linear isometry $\map hUV$.
In particular, $U$ and $V$ are linearly isometric.
\end{enumerate}
\end{prop}

\begin{pf}
Let $U$ be a Banach space of $\ud(\fD_{<\kappa})$ and fix Banach spaces $X\subs Y$ as in (0). Fix an isometric embedding $\map fXU$.
Choose a continuous chain $\sett{X_\al}{\al < \kappa}$ of closed subspaces of $Y$ so that $X_0 = X$, $X_\al \in \fD_{<\kappa}$ and $\bigcup_{\al < \kappa}X_\al$ is dense in $Y$.
Recall that a ``continuous chain" means that $X_\delta$ is the closure of $\bigcup_{\xi < \delta}X_\xi$ for every limit ordinal $\delta < \kappa$.
Using the definition of universal disposition, construct inductively a sequence of linear isometric embeddings $\map {f_\al}{X_\al}U$ so that $f_0 = f$ and $f_\beta \rest X_\al = f_\al$ whenever $\al < \beta$.
At limit steps we use the continuity of the chain.
The unique map $\map {f_\kappa}XU$ satisfying $f_\kappa \rest X_\al = f_\al$ for $\al < \kappa$ is an isometric embedding extending $f$. This shows both (0) and (1), since we may take $X=0$.

The proof of (2) is a standard back-and-forth argument.
Namely, let $\sett{U_\al}{\al < \kappa}$ and $\sett{V_\al}{\al < \kappa}$ be continuous chains of closed subspaces of $U$ and $V$ respectively, such that $U_\al$, $V_\al$ are of density $<\kappa$ for $\al < \kappa$ and $U = \cl (\bigcup_{\al < \kappa}U_\al )$, $V = \cl (\bigcup_{\al < \kappa}V_\al )$ (note that the closure is irrelevant if $\kappa > \aleph_0$).
Furthermore, we assume that $U_0 = X$ and $V_0 = Y$.
Construct inductively isometric embeddings $\map {f_\xi}{U_{\al(\xi)}}{V_{\beta(\xi)}}$ and $\map {g_\xi}{V_{\beta(\xi)}}{U_{\al(\xi+1)}}$ so that $f_0 = f$, $g_\xi \cmp f_\xi$ is the inclusion $U_{\al(\xi)} \subs U_{\al(\xi+1)}$, and $f_{\xi+1} \cmp g_\xi$ is the inclusion $V_{\beta(\xi)} \subs V_{\beta(\xi+1)}$ for each $\xi < \kappa$.
The limit steps make no trouble because of the continuity of both chains.
The regularity of $\kappa$ is used for the fact that every subspace of $U$ (or $V$, respectively) density $< \kappa$ is contained in some $U_\al$ (or $V_\beta$, respectively).
The ``limit" operators $\map {f_\kappa}UV$ and $\map {g_\kappa}VU$ are bijective linear isometries because $f_\kappa \cmp g_\kappa = \id V$ and $g_\kappa \cmp f_\kappa = \id U$.
Finally, note that $f_\kappa$ extends $f$, which completes the proof of (2).
\end{pf}

The next result is a special case of more general constructions, known in model theory (see, e.g., \jon~\cite{Jon}).
For Banach spaces this can be found in \cite{ACCGMud} and \cite{K_flims}.

\begin{tw}\label{Tgnoir}
Let $\mu$ be a cardinal and let $\kappa$ be an uncountable cardinal. 
Let $X$ be a Banach space of density $\loe\mu$.
Then there exists a Banach space $Y\sups X$ of density $\mu^{<\kappa}$ that is of universal disposition for spaces of density $<\kappa$.
\end{tw}

\begin{pf}
The space $Y$ will be constructed by using The Pushout Lemma.
So, we need to compute first, how many ``possibilities" we have.
The idea is that we first want to extend $X$ to a bigger Banach space $Z(X)$ such that every isometric embedding $\map fEF$ with $E\subs X$ and $F$ of density $<\kappa$ is \emph{realized} in $Z(X)$, that is, there exists an isometric embedding $\map gF{Z(X)}$ such that $g(f(x)) = x$ for $x\in E$.

Given an isometric embedding $\map fEF$ such that $E\subs X$, let $P(X,f)$ be the resulting Banach space of the pushout of $f$ and the inclusion $E\subs X$.
Clearly, the density of $P(X,f)$ is the maximum of $\dens X$ and $\dens F$.

Observe that there are at most $\mu^{<\kappa}$ closed subspaces of $X$ of density $<\kappa$.
This follows from the fact that the cardinality of $X$ is $\loe \mu^{\aleph_0}$.
Now, given two spaces $E$ and $F$ of density $\lam < \kappa$, the cardinality of the set of all isometric embeddings of $E$ into $F$ cannot exceed $\lam^\lam = 2^\lam \loe \mu^\lam$.
Finally, note that there are at most $2^{<\kappa}$ isometric types of Banach spaces of density $< \kappa$.
Here we use the fact that $\kappa$ is uncountable and therefore $2^{<\kappa}\goe \cont$.

It follows that there is a family $\Ef$ of cardinality $\loe\mu^{<\kappa}$ consisting of isometric embeddings $\map fEF$ with $E\subs X$, the density of $F$ is $<\kappa$ and every isometric embedding $\map gGH$ satisfying these conditions is isometric to some element of $\Ef$.
Write $\Ef = \sett{f_\xi}{\xi < \lam}$, where $\lam = |\Ef|$.
Construct inductively a continuous chain of Banach spaces $\sett{X_\xi}{\xi < \lam}$, starting with $X_0 = X$ and setting $X_{\xi+1} = P(X_\xi,f_\xi)$.
Let $Z(X) = X_{\lam}$, the completion of the union of $\sett{X_\al}{\al < \lam}$.

Note that every isometry from a subspace of $X$ of density $<\kappa$ into a space of density $<\kappa$ is realized in $Z(X)$, because we have taken care of all possibilities.
Furthermore, observe that for $\mu_1 = \dens{Z(X)}$ we have that $\mu_1^{<\kappa} = \mu^{<\kappa}$.
This follows from the fact that $\mu_1 \loe \mu^{<\kappa}$ and $(\mu^{<\kappa})^{<\kappa} = \mu^{<\kappa}$.

By the remark above, we can repeat this procedure up to $\mu^{<\kappa}$ many times, not enlarging the density.
That is, we construct a continuous chain of Banach spaces $\sett{Z_\al}{\al < \theta}$, where $\theta = \mu^{<\kappa}$, $Z_0 = X$ and $Z_{\al+1} = Z(Z_\al)$ for $\al < \theta$.
We claim that the resulting Banach space $Y = \bigcup_{\al < \theta}Z_\al$ is of universal disposition for spaces of density $<\kappa$.
Its density is exactly $\mu^{<\kappa}$.
The only thing is to check that the cofinality of $\theta$ is $\goe \kappa$.
In fact, a well known fact from cardinal arithmetic says that $\theta^{\cf(\theta)} > \theta$.
On the other hand, $\theta^\lam = \theta$ for every $\theta < \kappa$.
Thus, indeed, the cofinality of $\theta$ is $\goe \kappa$ and therefore every subspace of $Y$ that is of density $< \kappa$ is actually contained in some $Z_\al$.
This completes the proof.
\end{pf}

Since $\cont^{<\aleph_1} = \cont^{\aleph_0} = \cont$, we obtain the following corollary, without extra assumptions on cardinal arithmetic.

\begin{wn}[{\cite{ACCGMud}}]
There exists a Banach space of density $\cont$ that is of universal disposition for separable Banach spaces.
\end{wn}

The arguments from the last part of the proof of Theorem~\ref{Tgnoir} show that the construction could be somewhat optimized. Namely, since we know that $\mu^{<\kappa}$ has cofinality $\goe\kappa$ and clearly $\mu^{<\kappa} \goe 2^{<\kappa} \goe \kappa$, we conclude that either $\mu^{<\kappa} = \kappa$ and $\kappa$ is a regular cardinal, or else $\mu^{<\kappa} \goe \kappa^+$ and $\kappa^+$ is always a regular cardinal. Thus, the space $Y$ can be constructed as the union of a continuous chain of length either $\kappa$ (if $\kappa$ is regular) or $\kappa^+$ (if $\kappa$ is singular).
On the other hand, it is not clear whether taking the shorter chain we really obtain a different Banach space.

The theorem above does not say anything about uniqueness.
The only known fact, coming from the general \fra-\jon\ theory, is as follows.

\begin{tw}\label{Twolocog}
Let $\kappa$ be an uncountable cardinal satisfying $\kappa^{<\kappa} = \kappa$.
Then there exists a unique, up to isometry, Banach space $\V_\kappa$ of density $\kappa$ and of universal disposition for Banach spaces of density $<\kappa$.
Furthermore, every partial isometry between subspaces of $\V_\kappa$ of density $<\kappa$ extends to a bijective isometry of $\V_\kappa$.
\end{tw}

\begin{pf}
The existence of $\V_\kappa$ is an application of Theorem~\ref{Tgnoir} with $\mu = \kappa$.
The second statement and the uniqueness of $\V_\kappa$ follow from Proposition~\ref{Prbeou}(2).
\end{pf}

Note that Theorem~\ref{Tgnoir} shows the existence of strong \Gurarii\ spaces.
In fact, all spaces that are $\ud(\fD_{<\kappa})$ are strong \Gurarii, but on the other hand one can construct a strong \Gurarii\ space using pushouts with finite-dimensional spaces only.
As proved in \cite{ACCGMud}, such a space is not $\udsep$. We explain the details in Section~\ref{Sroleceezero}, showing that it is even not universal for spaces of density $\aleph_1$.

Note that the ``pushout construction" can be continued ``forever".
In other words, there is no upper bound for the density of a strong \Gurarii\ space.
In fact, a well-known property of infinite cardinals is that if $\mu^{<\kappa} = \mu$ then $(\mu^+)^{<\kappa} = \mu^+$, therefore one can use Theorem~\ref{Tgnoir} to construct spaces of $\ud(\fD_{<\kappa})$ that are of densities $\mu^+$, $\mu^{++}$, and so on. The problem of existence arises when one reaches a limit cardinal, however it can always be ``skipped", replaced by its successor.

In view of the recent results of Avil\'es and Brech \cite{AviBre}, a strong \Gurarii\ space of density $\cont$ constructed by pushouts is in some sense unique, as long as $\cont$ is a regular cardinal.

\section{On the structure of strong \Gurarii\ spaces}\label{Sergowrgwrg}

The following had already been observed by \Gurarii. The proof comes from his work~\cite{Gurarii}.

\begin{prop}
No separable Banach space can be a strong \Gurarii\ space.
\end{prop}

\begin{pf}
Suppose $U$ is a separable strong \Gurarii\ space.
For every two points $a$, $b$ on the unit sphere of $U$ there exists a unique linear isometry $\map f{X_a}{X_b}$ satisfying $f(a) = b$, where $X_a$, $X_b$ are linear spans of $\sn a$ and $\sn b$ respectively.
Applying Proposition~\ref{Prbeou}(2), we conclude that for every two points $a$, $b$ on the unit sphere of $U$ there exists a bijective isometry $h$ of $U$ such that $h(a) = b$.

Now, using a theorem of Mazur on the existence of smooth points on the unit sphere in every separable Banach space, we deduce that every point on the unit sphere of $U$ is smooth.
Recall that $p \in \usphere U$ is \emph{smooth} if there exists only one functional $\phi \in U^*$ such that $\norm \phi = 1 = \phi(p)$.

Finally, we get a contradiction by applying Proposition~\ref{Prbeou}(1) which says that every separable Banach space is isometric to a subspace of $U$; in particular the unit sphere of $U$ must contain non-smooth points.
Note that a point that is non-smooth in a subspace of $U$ cannot be smooth in $U$, by the Hahn-Banach extension theorem.
\end{pf}

A Banach space $X$ is called \emph{transitive} if for every $a,b$ in the unit sphere of $X$ there exists a bijective isometry $\map hXX$ such that $h(a)=b$.
The argument above shows that a transitive separable space must be smooth.
This is closely related to the Mazur's rotation problem:
Does there exist a separable transitive Banach space, different from the Hilbert space?
According to our knowledge, this problem is still open.

\begin{prop}\label{PLzipinus}
The \Gurarii\ space is not 1-injective for finite-dimen\-sio\-nal Banach spaces.
\end{prop}

\begin{pf}
According to \cite[Example 6.2]{ZippinHB}, there exists a Banach space $E = C(K)$, where $K$ is a metric compact space, that is not $1$-injective for finite-dimensional Banach spaces.
Every $C(K)$ space is a Lindenstrauss space (see \cite{MicPel}), therefore by Theorem~\ref{TWojtaszz} the space $E$ is $1$-complemented in the \Gurarii\ space $\G$.
Finally, if $\G$ were $1$-injective for finite-dimensional spaces, then so would be $E$, a contradiction.
\end{pf}

The following negative result is in contrast to Theorem~\ref{Tprispris}.

\begin{tw}\label{Tmejkveten}
Let $E$ be a non-separable strong \Gurarii\ space and let $\Gee$ be a skeleton in $E$.
Then there exists $G\in \Gee$ that is not $1$-complemented in $E$.
\end{tw}

\begin{pf}
Suppose $\Gee$ is a skeleton in $E$ such that each $G\in\Gee$ is $1$-complemented in $E$.
By Theorem~\ref{Tgcerugh} and Proposition~\ref{Pskelegeres},
we may assume that each member of $\Gee$ is linearly isometric to the \Gurarii\ space $\G$.
We now claim that $\G$ is $1$-injective for finite-dimensional spaces, which in view of Proposition~\ref{PLzipinus} is a contradiction.

Fix finite-dimensional spaces $X \subs Y$ and fix an operator $\map fXE$ with $\norm f \loe 1$.
By the Pushout Lemma, there are a finite-dimensional space $W$, an isometric embedding $\map j{\img fX}W$ and a linear operator $\map gYW$ such that $\norm g \loe 1$ and $g \rest X = j \cmp f$.
There exists $G\in \Gee$ such that $\img fX \subs G$.
Let $\map PEE$ be a projection such that $\norm P = 1$ and $\img PE = G$.
Using the fact that $E$ is a strong \Gurarii\ space, we find an isometric embedding $\map kWE$ such that $k \cmp j$ is the inclusion $\img fX \subs E$. The operator $P \cmp k \cmp g$ is an extension of $f$ and has norm $\loe1$.
\end{pf}

Note that exactly the same proof shows that $\G$ is not a strong \Gurarii\ space.  This argument does not use Mazur's theorem on the existence of smooth points.

Recall that a Banach space is \emph{weakly Lindel\"of determined} if its dual has a weak star continuous one-to-one linear operator into some $\Sigma$-product of the real lines, i.e., a linear topological space of the form
$$\Sigma(\Gamma) = \setof{x\in \Err^\Gamma}{|\setof{\gamma}{x(\gamma)\ne0}| \loe\aleph_0},$$
endowed with the product topology.
This class of Banach spaces contains all weakly compactly generated (in particular, all reflexive) spaces.
It is well known (see, e.g., \cite[Ch. 19]{KKLP}) that a weakly Lindel\"of determined Banach space always contains a skeleton of 1-complemented subspaces and this does not depend on the norm of the space (i.e. it holds after any renorming).
Thus, Theorem~\ref{Tmejkveten} gives the following

\begin{wn}
No strong \Gurarii\ space can be weakly Lindel\"of determined.
\end{wn}

One can go further and conclude that no strong \Gurarii\ space has a monotone (transfinite) Schauder basis (see, e.g., \cite{JordiS} for the definition and results on transfinite Schauder bases).
The reason is again that such a space has a skeleton of 1-complemented spaces (with standard monotone Schauder bases).
This property, however, is not preserved after renormings and indeed it is not clear whether there exists a strong \Gurarii\ space with any transfinite Schauder basis, or more generally, isomorphic to a space with a projectional resolution of the identity.

\section{The role of $c_0$}\label{Sroleceezero}

A well known theorem of Sobczyk~\cite{Sobczyk} says that $c_0$ is complemented in every separable Banach space.
More precisely, for every isometric embedding $\map i{c_0}X$ with $X$ separable, there exists a linear operator $\map TX{c_0}$ satisfying $T \cmp i = \id{c_0}$ and $\norm T \loe 2$ (see, e.g., the proof of Sobczyk's theorem in \cite[Thm. 17.2]{KKLP}).
We are going to prove the same for the class of ``pushout generated" Banach spaces that includes some strong \Gurarii\ spaces (see \cite{ACCGMud} or remarks after the proof of Theorem~\ref{Twolocog} above). As a consequence, we answer Problem 1 from \cite{ACCGMud}.

The next fact explains why complementability of $c_0$ forces the space not to be of universal disposition for Banach spaces of density $\loe\aleph_1$.
For this aim we need to know the fact that Sobczyk's theorem fails for Banach spaces of density $\aleph_1$ (regardless of the validity of the continuum hypothesis).

Recall that a family $\Aaa$ of infinite subsets of $\Nat$ is \emph{almost disjoint} if $A\cap B$ is finite for every $A\ne B$ in $\Aaa$.
There is a natural locally compact topology on $\Nat \cup \Aaa$ whose base consists of all the singletons of $\Nat$ and all sets of the form $\sn A \cup (A\setminus F)$ with $F\subs \Nat$ finite.
Let $K_\Aaa$ be the one-point compactification of this space.
In the literature, spaces of the form $K_\Aaa$ are often called \emph{Mr\'owka compacta}.
Notice that $C(K_\Aaa)$ has a natural isometric copy of $c_0$; the standard basis consists of all characteristic functions of the singletons of $\Nat$.
This copy of $c_0$ is not complemented in $C(K_\Aaa)$, unless $\Aaa$ is countable.
For the proof, see \cite[Thm. 17.3]{KKLP}.
Clearly, $\Aaa$ can be taken so that $|\Aaa| = \aleph_1$ and therefore $c_0$ is not complemented in a space of density $\aleph_1$.

\begin{prop}\label{Pcezerofj}
Let $X$ be a Banach space of $\udsep$. Then no copy of $c_0$ can be complemented in $X$.
\end{prop}

\begin{pf}
Let $c_0 \subs Z$, where $Z = C(K_\Aaa)$ for some almost disjoint family $\Aaa$ of cardinality $\aleph_1$.
Let $E \subs X$ be isomorphic to $c_0$ and let $\map f{c_0}E$ be an isomorphism.
Using Lemma~\ref{Lertnort}, find an equivalent norm on $Z$ such that $f$ becomes an isometry.
By Proposition~\ref{Prbeou}(0), there is an isometry $\map gZX$ such that $g\rest c_0 = f$.
It is now clear that $E$ cannot be complemented in $\img gZ \subs X$, therefore it cannot be complemented in $X$.
\end{pf}

We are now going to show that Sobczyk's theorem holds in a class of Banach spaces containing strong \Gurarii\ spaces of arbitrarily large density.

\begin{df}\label{Dpngoggb}
Let $\pogfd$ denote the class of all Banach spaces that can be obtained as the limit (i.e. the completion of the union) of a transfinite chain $\sett{X_\al}{\al<\rho}$ such that $X_0$ is separable, $X_\delta = \cl(\bigcup_{\xi<\delta}X_\xi)$ for every limit ordinal $\delta<\rho$ and for each $\al<\rho$, the space $X_{\al+1}$ comes from the pushout square
$$\xymatrix{
X_\al \ar[r]^\subs & X_{\al+1} \\
E_\al \ar[u]^{j_\al}\ar[r]_\subs & F_\al \ar[u]
}$$
where $E_\al \subs F_\al$ are finite-dimensional spaces and $j_\al$ is an isometric embedding.
More specifically, we shall write $X\in \pogfd(Y)$ if $X$ is the limit of a chain as above, in which $Y=X_0$.
\end{df}

As mentioned before, it has been proved in \cite{ACCGMud} that the class $\pogfd$ contains strong \Gurarii\ spaces.

Before proving our result, we need the following lemma, which can be easily deduced from a variant of \cite[Lemma 20]{AviBre} involving finite-dimensional spaces.

\begin{lm}\label{Labeorgrg}
Let $Z$ be a separable subspace of a space $X\in\pogfd$. Then there exists a separable space $Y\subs X$ such that $Z \subs Y$ and $X \in\pogfd(Y)$.
\end{lm}

\begin{tw}
Let $X \in \pogfd$. Then every copy of $c_0$ is complemented in $X$.
\end{tw}

\begin{pf}
Let $C\subs X$ be isometric to $c_0$. By Lemma~\ref{Labeorgrg}, we may assume that $C\subs X_0$ for some separable space $X_0$ such that $X = \cl(\bigcup_{\xi<\rho}X_\xi)$, where the chain $\sett{X_\xi}{\xi<\rho}$ satisfies the conditions in Definition~\ref{Dpngoggb}.
By Sobczyk's theorem, there exists a projection $\map P{X_0}C$ with $\norm P \loe2$.

Set $P_0 = P$. We now construct inductively projections $\map {P_\al}{X_\al}C$ so that $P_\beta$ extends $P_\al$ whenever $\beta > \al$ and $\norm {P_\al} = \norm P$ for every $\al$.
Suppose $P_\xi$ have been constructed for $\xi < \al$.
If $\al$ is a limit ordinal, we define $P_\al$ to be the pointwise limit of $\sett{P_{\xi_n}}{\ntr}$, where $\xi_0 < \xi_1 < \dots < \al$ converges to $\al$.
Here we have used the fact that $X_\al$ is the closure of $\bigcup_{\ntr}X_{\xi_n}$.

Now suppose $\al = \eta + 1$ and fix a pushout square
$$\xymatrix{
X_\eta \ar[r]^\subs & X_{\al} \\
E \ar[u]^{j}\ar[r]_\subs & F \ar[u]_k
}$$
defining $X_{\al}$, with finite-dimensional spaces $E,F$.
Using the fact that $c_0$ is $1$-injective for finite-dimensional spaces, we find a linear operator $\map TFC$ satisfying $T\rest E = P_\eta \cmp j$ and $\norm T = \norm {P_\eta \cmp j} = \norm {P_\eta}$.
By the pushout property, there exists a unique operator $\map{P_\al}{X_\al}C$ satisfying $P_\al \rest X_\eta = P_\eta$, $P_\al \cmp k = T$ and $\norm {P_\al} = \norm {P_\eta}$.

Finally, $P = \lim_{\xi < \rho}P_\xi$ is the required projection.
\end{pf}

It has been shown in \cite{ACCGMud} (with almost the same arguments) that if $X\in\pogfd(Y)$, where $Y$ is linearly isometric to $c_0$, then $Y$ is $1$-complemented in $X$.

\begin{wn}
Let $X\in\pogfd$. Then $X$ cannot contain any isomorphic copy of $C(K_\Aaa)$, where $\Aaa$ is an almost disjoint family of infinite subsets $\Nat$ and $|\Aaa| = \aleph_1$.
\end{wn}

This answers Problem 1 from \cite{ACCGMud}: There exist strong \Gurarii\ spaces (of arbitrarily large density) that are not universal for Banach spaces of density $\aleph_1$.

\section{Final remarks and open problems}\label{SecTheLastOne}

Below we collect some open questions; some of them are motivated by the results described in previous sections.

\paragraph{Minimal density.}
It is not clear what the minimal density of a strong \Gurarii\ space is.
The only known bound is the continuum. A more concrete question is:

\begin{question}
Does there exist, without extra set-theoretic assumptions, a strong \Gurarii\ space of density $\aleph_1$?
\end{question}

\begin{question}
Assuming $\cont < \aleph_\omega$, does there exist a strong \Gurarii\ space of density $\aleph_\omega$?
\end{question}

Note that $\aleph_\omega$ is the smallest singular cardinal and it has cofinality $\omega$; therefore always $\cont \ne \aleph_\omega$.

\paragraph{Schauder bases.}
A Banach space with a PRI and of density $\aleph_1$ has a countably 1-norming Markushevich basis (see, e.g., \cite[Section 17.8]{KKLP}). A Markushevich basis is a natural ``non-separable" generalization of Schauder bases, yet it exists in every separable Banach space. Theorem~\ref{Tprispris} motivates the following

\begin{question}
Does there exist a \Gurarii\ space of density $\aleph_1$ with a monotone transfinite Schauder basis?
\end{question}

Note that by Theorem~\ref{Tmejkveten} such a space cannot be strong \Gurarii.
Let us mention that some of the ``generic" Banach spaces constructed in \cite{JordiS} are \Gurarii, although none of them has a transfinite Schauder basis.

\begin{question}
Does there exist a strong \Gurarii\ space, isomorphic to a Banach space with a PRI?
\end{question}

\begin{question}
Does there exist a weakly Lindel\"of determined (or even better: weakly compactly generated) \Gurarii\ space?
\end{question}

Again, this cannot be a strong \Gurarii\ space.
Note that every weakly Lindel\"of determined Banach space has a countably 1-norming Markushevich basis.

\paragraph{Renormings.}
Recall that a norm $\anorm$ is \emph{rotund} if 
$\norm{x+y} = 2\norm x = 2\norm y$ implies $x = y$.
A \emph{rotund renorming} is an equivalent norm that is rotund.
Many non-separable Banach spaces have rotund renormings, for a general treatment we refer to the book \cite{DGZ}.
A result of Zizler~\cite{Ziz} says that the existence of a renorming stronger than rotund (namely: locally uniformly rotund) is preserved by a PRI.
In particular, every Banach space of density $\aleph_1$ and with a PRI has a rotund renorming.
In view of Theorem~\ref{Tprispris}, there exist non-separable \Gurarii\ spaces admitting a rotund renorming.
This suggests:

\begin{question}
Does there exist a strong \Gurarii\ space with a rotund renorming?
\end{question}

A typical example of a Banach space with no rotund renorming is $\ela \by c_0$ (see \cite{DGZ}).
Unfortunately, this space has density $\cont$ and the following interesting question, due to Christina Brech, seems to be open.

\begin{question}\label{Qunwgo}
Does there exist, without extra set-theoretic assumptions, a Banach space $X$ of density exactly $\aleph_1$ and with no rotund renorming?
\end{question}

A positive answer to this question would yield a simple and direct proof of the following result.

\begin{tw}\label{Tfgrwgr}
No Banach space of universal disposition for separable spaces can have a rotund renorming.
\end{tw}

Indeed, a space of $\udsep$ contains copies of all Banach spaces of density $\aleph_1$, so all of them would have to admit rotund renormings.
Assuming CH, this gives a contradiction.
Still, the statement above is a theorem.
For readers familiar with the technique of forcing, we sketch a ``metamathematical" proof, involving absoluteness.

\begin{pf}
Suppose the statement above is not a theorem, i.e. it is not a consequence of the usual axioms of set theory.
By G\"odel's completeness, there exists a model of set theory $\V$ that contains a Banach space $X$ of $\udsep$ with rotund renorming.
There exists an extension $\W$ of $\V$ (obtained by forcing) such that $\W$ is a model of set theory in which the continuum hypothesis holds and, moreover, for every function $\map \phi\omega S$ in $\W$ if $S \in \V$ then $\phi\in \V$.
The last property of $\W$ implies that $X$ is a Banach space in $\W$ and it is of $\udsep$.
The latter fact is because $\W$ does not contain ``new" separable Banach spaces.
Finally, $X$ still has a rotund renorming, since this property is preserved.
This leads to a contradiction, since in $\W$ the space $X$ contains a copy of $\ela \by c_0$.
\end{pf}

\paragraph{Ultra-homogeneity.}
Let $\fK$ be a class of Banach spaces.
We say that a Banach space $X$ is \emph{homogeneous} with respect to $\fK$ if every bijective isometry between two subspaces of $X$ that are in class $\fK$ extends to an isometry of $X$ onto itself.
If $\fK$ contains all 1-dimensional subspaces of $X$, homogeneity implies transitivity.
In fact, the difficulty of Mazur's problem on rotations exhibits the fact that 
so far the Hilbert space is the only known example of a separable Banach space homogeneous for finite-dimensional spaces.
Now let $\fK = \fD_{<\kappa}$, the class of all Banach spaces of density $<\kappa$.
Proposition~\ref{Prbeou}(2) says that every space of $\ud(\fK)$ and of density $\kappa$ is homogeneous with respect to $\fK$.
On the other hand, in view of the results of \cite{AviBre}, there exist (arbitrarily large) Banach spaces of $\udsep$ that are homogeneous with respect to separable subspaces.
It is not clear what happens with strong \Gurarii\ spaces.
\begin{question}
Does there exist a strong \Gurarii\ space, homogeneous with respect to finite-dimensional subspaces and not of universal disposition for separable spaces?
\end{question}

\begin{question}
Does there exist a strong \Gurarii\ space that is not homogeneous for finite-dimensional spaces?
\end{question}

In fact, we do not know the answer to a more general question:

\begin{question}
Does there exist a Banach space of $\ud(\fD_{<\kappa})$ that is not homogeneous with respect to $\fD_{<\kappa}$?
\end{question}

We finish with the following problem whose solution may lead to a better understanding of Mazur's rotation problem.

\begin{problem}
Find a class $\fK$ of finite-dimensional Banach spaces with the following properties:
\begin{enumerate}
	\item[(i)] $\fK$ is hereditary (i.e. $X\subs Y\in \fK$ implies $X\in \fK$).
	\item[(ii)] All spaces in $\fK$ are smooth.
	\item[(iii)] For each $\Ntr$, $\fK$ contains a space of dimension $n$.
	\item[(iv)] $\fK$ has the amalgamation property. That is, given isometric embeddings $\map iZX$, $\map jZY$ with $X,Y\in \fK$, there exist $W\in \fK$ and isometric embeddings $\map {i'}XW$, $\map {j'}YW$ satisfying $j' \cmp j = i' \cmp i$.
	\item[(v)] $\fK$ is not dense (with respect to the Banach-Mazur distance) in the class of all finite-dimensional Banach spaces.
	\item[(vi)] $\fK$ is not the class of Euclidean spaces.
\end{enumerate}
\end{problem}

Actually, it is desirable to replace condition (vi) by a formally stronger one: $\fK$ contains a chain $\ciag X$ such that the completion of $\bigcup_{\ntr}X_n$ is not isomorphic to the Hilbert space.

Having such a class $\fK$, one would be able to construct a Banach space $\G_\fK$ satisfying the definition of the \Gurarii\ space for finite-dimensional spaces from class $\fK$ only.
If the class $\fK$ had an additional property that $\G_\fK$ remains smooth (which does not follow from condition (ii)), \Gurarii's argument would not be applicable for showing that $\G_\fK$ is not transitive.
In any case, $\G_\fK$ would be a new Banach space ``almost" homogeneous with respect to its finite-dimensional subspaces and not isomorphic to the Hilbert space.

\separator

\paragraph{Acknowledgments} We would like to thank Antonio Avil\'es and Jes\'us Castillo for useful remarks and comments.

\end{document}